\documentclass[12pt]{report}

\usepackage{amsfonts,amssymb,amsmath,amscd,latexsym,epsfig,graphicx}

\begin{document}

\begin{flushright}
3 may 2005 
\end{flushright}

 \begin{center}
 \textbf{\Large On Lie bialgebras of loops on orientable surfaces }
 \end{center}

\begin{center}
\textbf {Attilio Le Donne}\footnote{Tel.: +39 06 49913249; fax:  +39 06 44701007 

E-mail address: ledonne@mat.uniroma1.it} 
\end{center}
\vskip0.5cm
Dipartimento di Matematica, Universit\' a di Roma "La Sapienza"  , Piazzale Aldo Moro 1, I-00185, Roma, Italia
\vskip0.5cm

\begin{flushleft}
 \textbf{Abstract}
\end{flushleft}                       
Goldman (Invent. Math. 85(2) (1986) 263) and Turaev (Ann. Sci. Ecole Norm. Sup. (4) 24 (6)(1991) 635) found a Lie bialgebra structure on the vector space generated by non-trivial free homotopy classes of loops on an orientable surface.
Chas ( Combinatorial Lie bialgebras of curves on surfaces, Topology 43 (2004) 543), by the aid of the computer, found  a negative answer  to  Turaev's question about the characterization of multiples of simple classes in terms  of the cobracket, in every surface of   negative Euler characteristic and   positive genus. However, she  left open  Turaev's conjecture, namely  if,  for genus zero,  every class  with cobracket zero is a  multiple of a simple class.
The aim of this paper is to  give a positive answer  to this conjecture.
\vskip0,5cm
MSC: primary 57M99, secondary 17B62.
\vskip0,5cm
Keywords: Surfaces; Homotopy  classes; Lie bialgebras.

\vskip0,5cm
\begin{flushleft}
 \textbf{1. Introduction.}
\end{flushleft}   
                                             
In 1989, Turaev [3] introduced for an oriented surface $\Sigma$, a Lie cobracket in the module $Z(\Sigma)$ generated by the homotopy classes of non-trivial loops on $\Sigma$. 

The cobracket complements the Goldman Lie bracket [2] in $Z(\Sigma)$ and makes it a Lie bialgebra in the sense of Drinfeld.

Every multiple of a simple class has cobracket zero, so  Turaev [4] asked whether it is true  that if the cobracket is zero then the class is a multiple of a simple class.

 Chas [1], giving a combinatorial presentations of the  Lie bialgebra by reduced cyclic words, found, by the aid of the computer, that in every surface of     negative Euler characteristic and   positive genus there exist classes with cobracket zero which are not multiples of simple classes. 

Always with the aid of the computer, she checked that Turaev's characterization holds for
words of length up to sixteen on the pair of pants (i.e. a two sphere with three punctures) so that she suggested that it is  possible for genus
zero surfaces.

Here we answer affirmatively this conjecture by showing that:

\begin{flushleft}
\emph{In a genus zero surface every class with cobracket zero is a multiple of a simple class.}
\end{flushleft}

\vskip0,5cm
\begin{flushleft}
\textbf{2. Curves and loops}
\end{flushleft}

Let $\Sigma$ be a smooth  oriented surface, possibly with boundary.

The circle $S^1$ is oriented counterclockwise, this will make all curves oriented.

Given  two  points $\alpha,\beta $ of the circle $S^1$ , we denote with $\alpha\beta $  (resp.$-\alpha\beta $) the  oriented embedded arc in $ S^1$ which starts at $\alpha$ (resp. $\beta $), going  in the positive, i.e.,    counterclockwise (resp. negative, i.e.,   clockwise)  direction and terminates at $ \beta $ (resp. $\alpha $), so in $\alpha\beta$ (resp.$-\alpha\beta $) we consider this order. If  $\alpha=\beta $,  $\alpha\beta $  (resp. $-\alpha\beta $) will be all $ S^1$ (resp.  $ -S^1$)  with counterclockwise (resp. clockwise)  direction.
Note that $-\alpha\beta $ is different from $\beta\alpha $, even if both start at  $\beta $ and end at $\alpha$.

 A    $curve$   in $\Sigma$ will be a piecewise smooth map $f$  from an arc $\alpha\beta $  to $ \Sigma$. If $f(\alpha) = f(\beta) $   the curve  is a $loop$.
 
Two homotopic curves $f$ and $g$ are said also $equivalent$, and we write $f \cong g$. A loop is trivial if it is null-homotopic.
 
From now on, $f$ will be a curve.

If the arc $ab $  is a subset of the domain  of $f$, then  $f_{ab}$ will be the restriction of  $f$ to  $ab $.

For two curves $f_{\alpha\beta}$ and $g_{\gamma\delta}$, so that $f(\beta)=g(\gamma)$, $f_{\alpha\beta}\cdot g_{\gamma\delta}$ is the usual composition of curves.

$f_{-\alpha\beta}$ will be the inverse curve of $f_{\alpha\beta}$. ( It start at  $f(\beta)$ and end at $f(\alpha)$).

For two equivalent curves  $f$ and $g$ we put  $g \prec f$ if  $ im (g) \subset im (f) $. (Note that  we can have $g \prec f$ also if $f$ and $g$ have not the same domain).

Denote with $Class (f)$ (resp. $RClass (f)$ ) the class of all $g$ homotopic to $ f$ (resp. $g$ such that $g\prec f$).

\vskip1cm
\begin{flushleft}
\textbf{3. Double and crossing points }
\end{flushleft}
 The point $p$ is a double point if $f^{-1} (p) = \{ x,y\},$ with $\ x\neq y  $; in this case  we put $x\approx y $.
 
If $f(a)$ is a double point, we call $a$ a double point too. 

We consider only curves with a finite number of double points and no other multiple point.
 
 Moving the curve in $V$ will mean  use a homotopy that  leaves  all point not in $V$ fixed.  

A double point $p$ is a  $ crossing $ point if,  for some neighborhood $V$ of $p$, we cannot eliminate the double point moving the curve   in $V$.

We denote $Dbl(f)$ (resp. $Cross(f)$ )  the set of double (resp. crossing ) points of $f$, and $\sharp Dbl(f)$ and $ \sharp Cross(f)$  their numbers.

An  $f$ is $generic$ if it is  smooth and has  only transversal ( i.e.,   with transversal tangents) double points; (hence  each double point is crossing, i.e.,   $\sharp Dbl(f)=\sharp Cross(f)$).
 
 We say  that  $ f$  is  $ minimal$ if it has  the minimal number of crossing points in $RClass (f)$, ( i.e.,    $\sharp Cross(f)\leq \sharp Cross(g)$ for each $g$ such that $g\prec f$).

\vskip0,3cm

The $Feynmann$  $diagram$ of $f_{\alpha\beta}$  is formed by the arc $\alpha\beta $  and all the straight segments $e _p=e_{ab}$, if $a<b$ in $\alpha\beta $, having endpoints the two points $a,b$ of $f^{-1}(p)$ , for each double point $p $ of $f$. We say that two different double points $p, q$ are linked if $e _p \cap e _q  \neq \emptyset$. 

We say that a Feynmann diagram is linked if there are two  double points which are  linked.

\vskip0,5cm
\begin{flushleft}
\textbf{4. Monogons, bigons and trigons }
\end{flushleft}

A straight segment $e_{ab}$ of the Feynmann diagram of $f$  forms a $monogon$ if  $f_{ab}$ is a null-homotopic loop.
\vskip0,1cm

Two linked straight segments $e_{ab}$ and $e_{cd}$ of the Feynmann diagram, with $a<c<b<d$, form a
 $direct \  bigon$ of  $f $ if   $f_{ac}\cong f_{bd}$. 
\vskip0,1cm

Two  non-linked straight segments $e_{ab}$ and $e_{cd}$ of the Feynmann diagram, with $a<c<d<b$, form a
 $inverse \  bigon$ of  $f $   if   $f_{ac} \cong f_{-bd}$.
\vskip0,1cm

A direct  bigon is a particular kind of a trigon.

Two linked straight segments $e_{ab}$ and $e_{cd}$ of the Feynmann diagram , with $a<c<b<d$, form a
 $trigon$ of  $f $ if   $  f_{bd}\cdot f_{cb}\cdot f_{ac}  \cong    f_{ac}   \cdot f_{cb}\cdot f_{bd} (=  f_{ad} )  $.

\vskip0,5cm
\begin{flushleft}
 \textbf{Lemma 1.}
\end{flushleft}

\textsl{Every  generic minimal non-trivial $f$ has neither  monogons nor bigons or trigons.}
\vskip0,5cm

\noindent \textbf{Proof.\ }Let $f=f_{\alpha\beta} $.

If $e_{ab}$ is a monogon of $f $ then put $g=f_{\alpha a}\cdot f_{ b \beta}$.
  
Let  $\alpha \leq a<b<c<d \leq \beta $.

If  $e_{ab}$ and $e_{cd}$ form an inverse  bigon of           $f $ then put $g=f_{\alpha a}\cdot f_{-cd}\cdot f_{bc}\cdot f_{-ab}\cdot f_{d \beta}$.

If $e_{ab}$ and $e_{cd}$  form direct bigon or a trigon of   $f $ then put $g=f_{\alpha a}\cdot f_{cd}\cdot f_{bc}\cdot f_{ab}\cdot f_{d  \beta}$.
\vskip0,1cm

In all cases, $g\prec f$ but $\sharp Cross(g)< \sharp Cross(f)$. Contradiction. $\diamond$

 \vskip0.5cm                                              

Let $\psi$ be a simple loop in $\mathbb{R}^2 $,   $D_\psi$ the bounded component of $\mathbb{R}^2 \setminus \psi$,  and $p_1, \ldots, p_n   \in D$. ( We say that  a curve  is inside $\psi$ if his image is in $cl(D_\psi)$ ).

Put $\Omega (\psi)= cl (D_\psi ) \setminus \  \left\{p_1, \ldots, p_n \right\}=  ( D_\psi \cup \psi )  \setminus \  \left\{p_1, \ldots, p_n  \right\}  $ 

\vskip0,3cm
\begin{flushleft}
 \textbf{Lemma 2.}
\end{flushleft}

\emph{Let  $ f  : \alpha\beta \longrightarrow \Omega(\psi) $   be a generic minimal curve  with $m$ double points so that:
\begin{enumerate}
	\item $ f( \alpha ) \in im (\psi) $,
\item  $\alpha$  and $ \beta$ are not double points of $f$, 
\item Each simple loop $f_{ab}$ is homotopic  to $\psi$ in $\Omega(\psi)$.
\end{enumerate}
Then there are $a_i,b_i\in\alpha\beta$, for each $i\leq m$, such that: 
\begin{itemize}
	\item  $\alpha<a_m<\ldots<a_1<b_1<\ldots<b_m<\beta$, and $a_i\approx b_i$, for $i\leq m$, (  i.e.,   $f$ has  a non linked Feynman diagram ),
	\item $\sigma_{i+1}= f_{a_ia_{i+1}}\cdot f_{b_{i+1}b_i }$, for $i< m$, and $\sigma_1= f_{a_1b_1}$ are simple equivalent loops (and hence homotopic to $\psi$),
	\item  $\sigma_i$ is inside $\sigma_{i+1}$,  $(\ im(\sigma_i) \cap im( \sigma_{i+1}) =\{f(a_i)\}= \{f( b_i)\}\  )$, for $i< m$;  $f(\beta)$  is outside  $\sigma_m$.
\end{itemize}
Hence $ f_{a_i b_i}\cong \psi^i$, for $i\leq m$}  ( see fig. 1)

\begin{center}

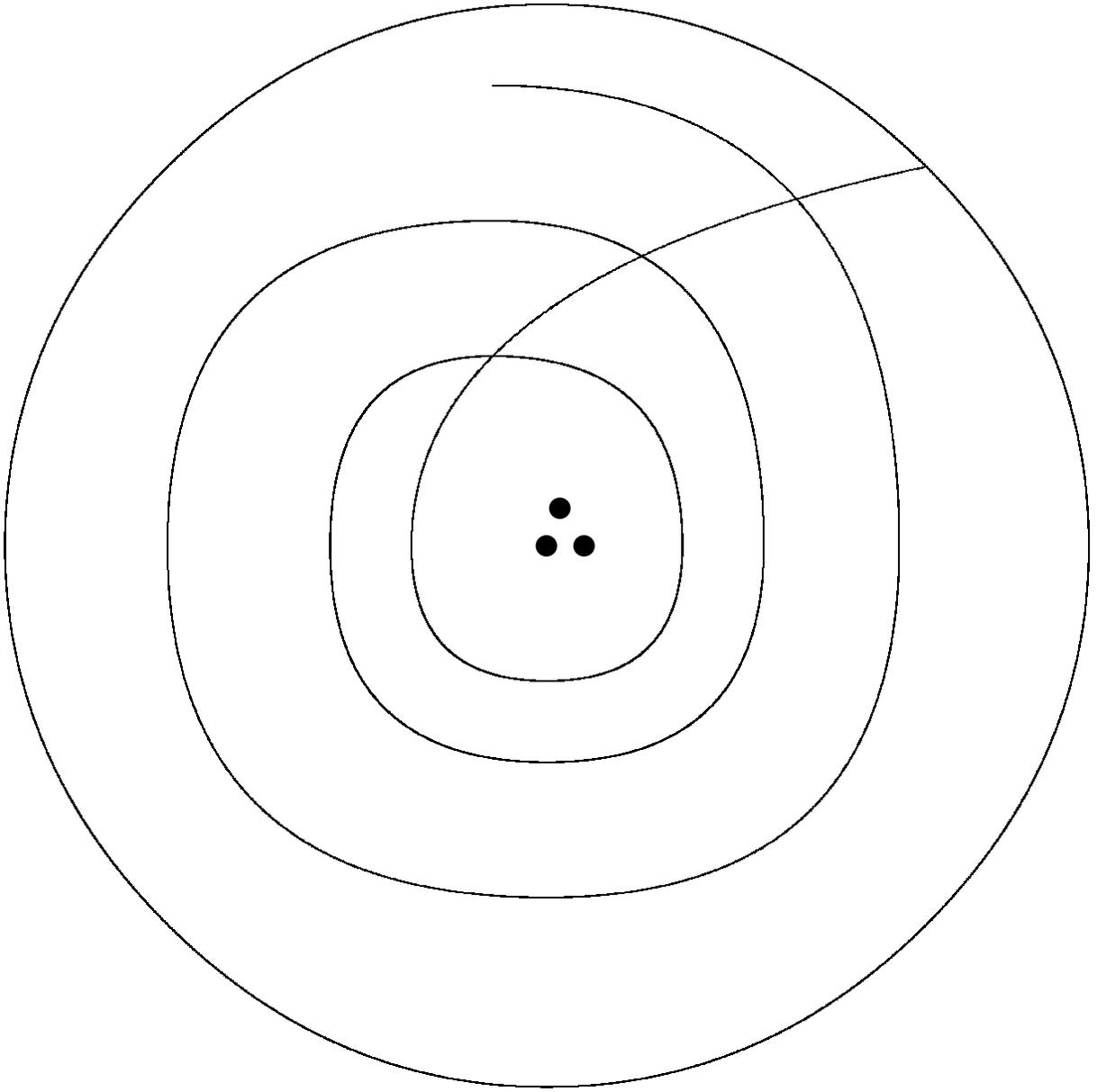
\textbf{Figure 1}
\end{center}
\begin{center}

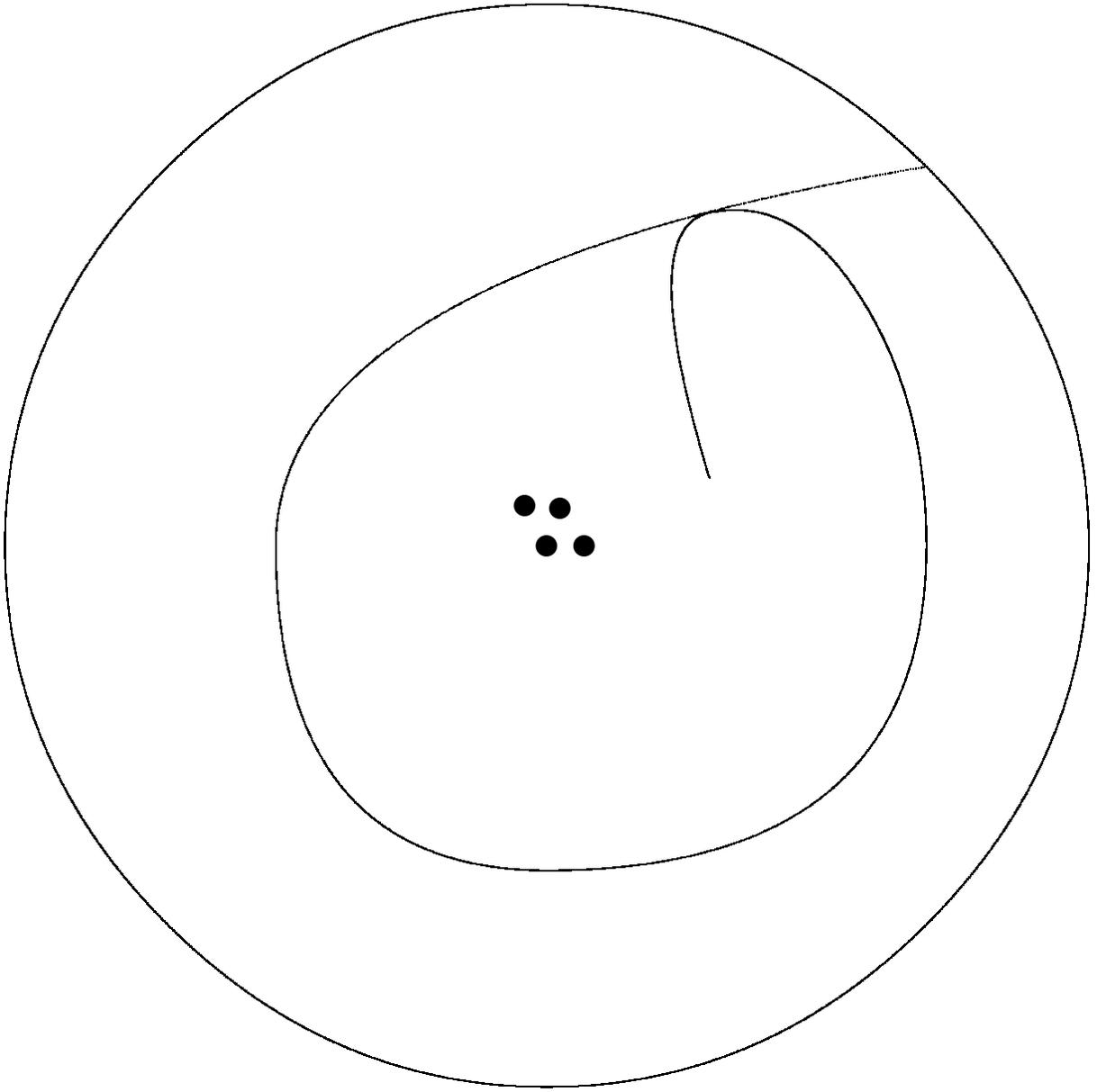
\textbf{Figure 2}
\end{center}

\begin{flushleft}
\textbf{Proof}
\end{flushleft}

The proof is by induction on $m$.

If  $m=1$ the only thing to demonstrate is that $f(\beta)$  is outside  $f_{a_1b_1}$. But, if it were not so, the double point $f(a_1)=f(b_1)$ would  not be a crossing point. (see fig. 2)

Let us assume  the lemma true for $m $ and  $f$,  with $ m+1$ crossing points, verifying  the hypothesis of lemma.
 
 Choose $\delta  \notin  Dbl(f) $ so that in $\delta \beta $ there is  only one  double point  that we call $b_{m+1}$ and let $h=f_{\alpha\delta}$.
 
 Clearly  $\sharp Dbl(h)= m$ and $h$ verifies the three properties of lemma and so there are $a_i,b_i\in\alpha\delta$ for each $i\leq m$, such that

\begin{itemize}
	\item $\alpha<a_m< \ldots<a_1<b_1<\ldots<b_m<\delta$, and  $a_i\approx b_i$, for $i\leq m$,
 	\item $\sigma_{i+1}= f_{a_i a_{i+1}}\cdot f_{b_{i+1} b_i }$, for $i< m$, and $\sigma_1= f_{a_1 b_1}$ are simple equivalent loops, 
 	\item  $\sigma_i$ is inside $\sigma_{i+1}$, for $i< m$;  $f(\delta)$  is outside  $\sigma_m$.
\end{itemize} 

Since in in $\delta \beta $ there is the only  double point  $b_{m+1}$ of $f$, we have that $f_{\delta \beta}$ cannot meet  $\sigma_i$ for  $i< m$ and so $f(\beta)$ is outside $\sigma_m$ or between  $\sigma_m$ and $\sigma_{m-1}$.

Let $a_{m+1}$ be such that  $a_{m+1}\approx b_{m+1}$.

So we have one of the following cases:
\begin{enumerate}
	\item $a_{m+1}\in \alpha  a_m  $ 
\item  $a_{m+1}\in a_m a_{m-1} $ \qquad ($a_2\in a_1b_1 $,\quad  if \ $m=1 $)
\item $a_{m+1}\in b_{m-1}b_m $ \qquad ( case not considered if \  $m=1 $)
\item $a_{m+1}\in b_m\delta$
\end{enumerate}

\begin{center}

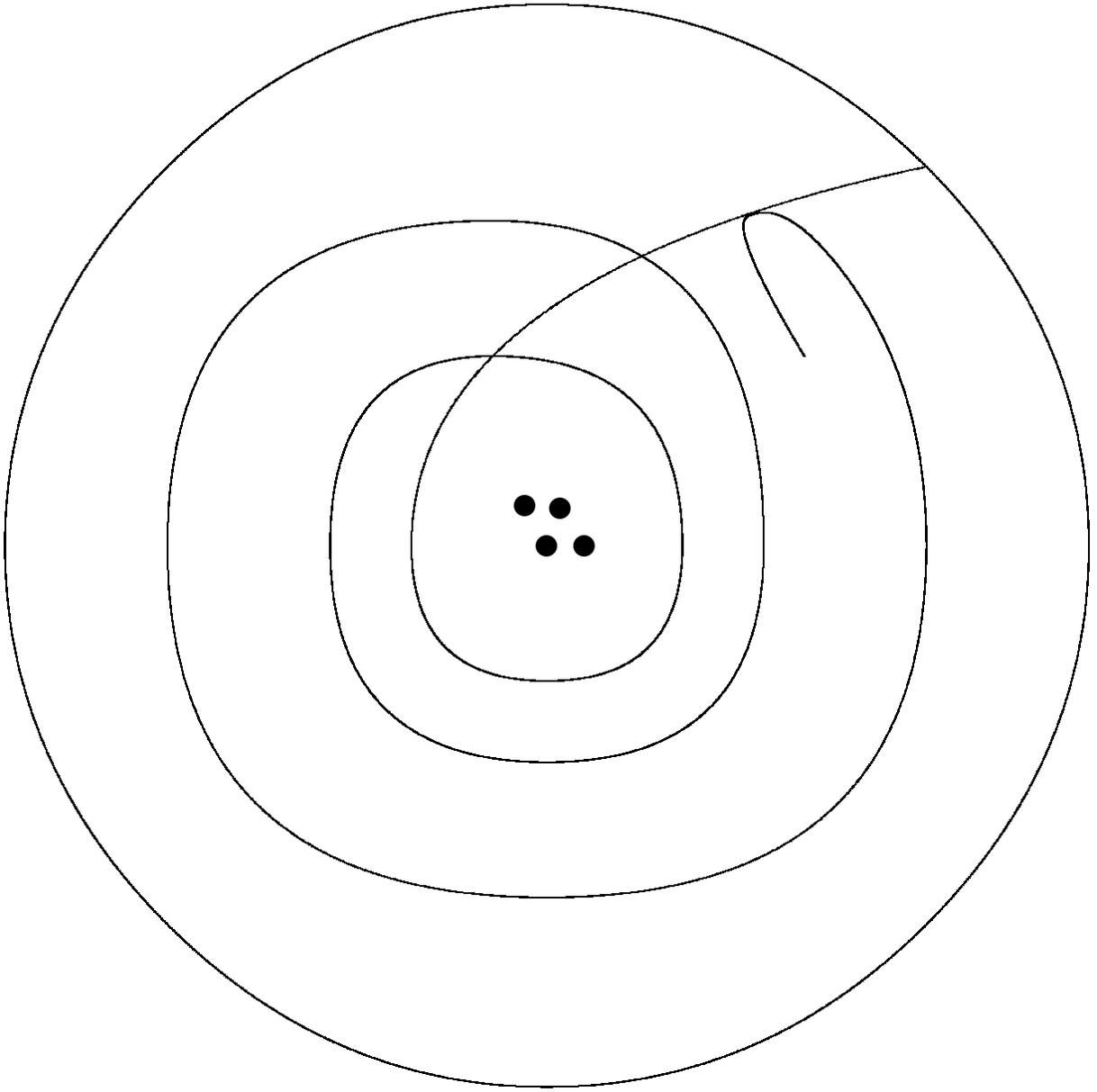
\textbf{Figure 3}
\end{center}

\begin{center}

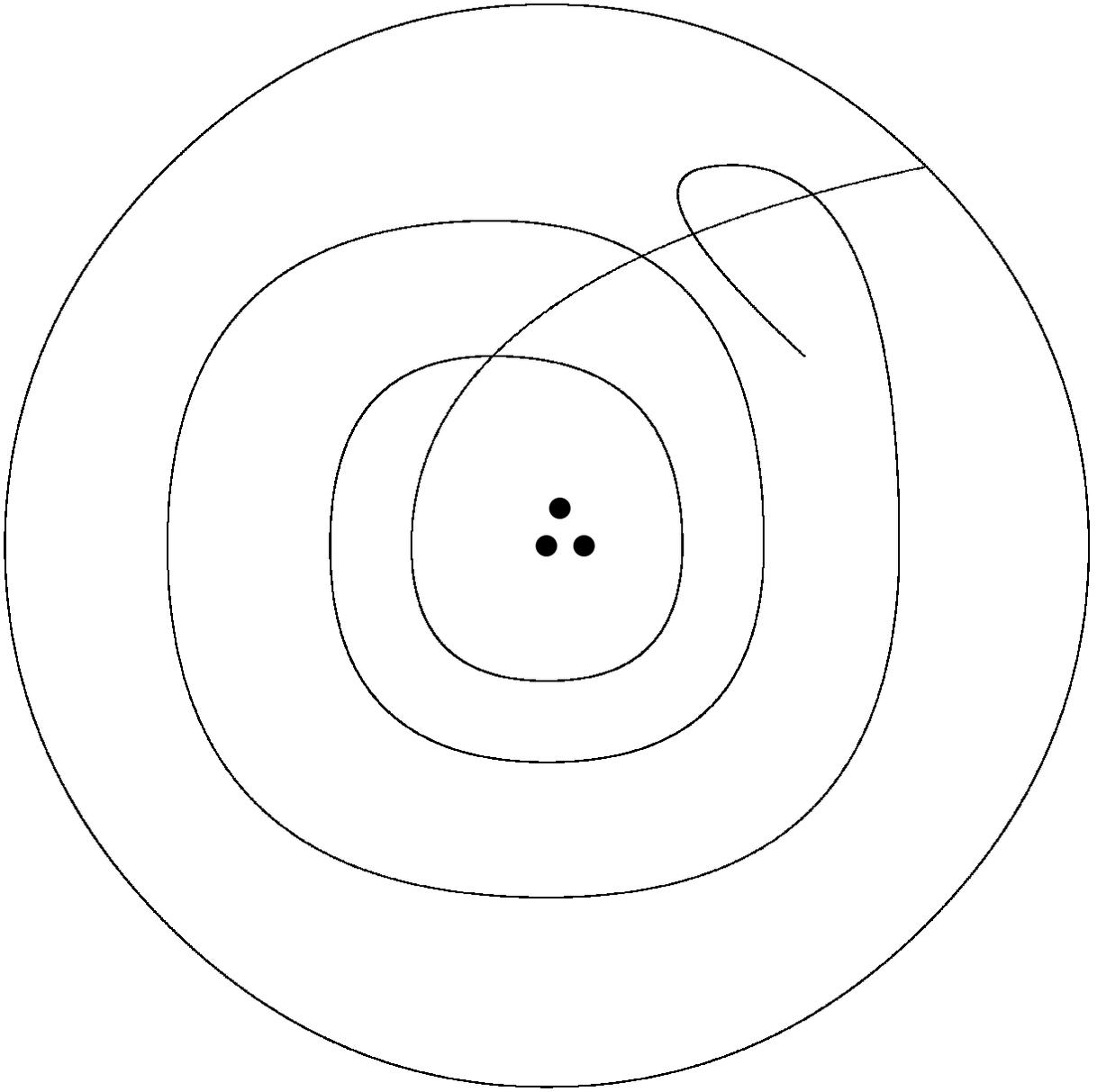
\textbf{Figure 4}
\end{center}

\begin{center}
\textbf{Figure 5}
\end{center}
      
\begin{center}
\textbf{Figure 6}
\end{center}

Put $e_i = e_{a_i b_i}$,  for $i\leq m+1$.

In case 2), $e_m $ and  $e_{m+1}$ would form a direct bigon, (see fig. 4); in fact $f_{a_m a_{m+1}}\cong f_{-b_m b_{m+1}}$.

In case 3),  $e_m $ and  $e_{m+1}$ would form a trigon, (see fig. 5), because 
$  f_{a_m a_{m+1}}\cdot f_{ a_{m+1}b_m} \cdot f_{b_mb_{m+1}}\cong 
f_{b_m b_{m+1}}\cdot f_{ a_{m+1}b_m} \cdot f_{a_m a_{m+1}}$;
 in fact  
 $(f_{a_m a_{m+1}}\cdot f_{ a_{m+1}b_m}) \cdot (f_{b_mb_{m+1}}\cdot f_{ a_{m+1}b_m})
  \cong (f_{b_m b_{m+1}}\cdot f_{ a_{m+1}b_m}) \cdot (f_{a_m a_{m+1}}\cdot f_{ a_{m+1}b_m})$, since their images are between $\psi$ and $\sigma_1$ and the fundamental group of  a circular crown  is commutative.
 
In case 4),   $e_{m+1}$ would be  a monogon, (see fig. 6).

In case 1), if $ \sigma_{m+1}= f_{ a_{m+1}a_m}  \cdot  f_{b_m b_{m+1}}$   were null homotopic then we would have an inverse bigon; otherwise       $ \sigma_{m+1}   \cong    \psi $.

It is still to be proved that   $f(\beta)$  is  outside  $\sigma_m$. But,
 as in case $m=1$, if   $f(\beta)$  were inside  $\sigma_m$, then the double point $f(a_{m+1})=f(b_{m+1})$ would not be  a crossing point, (see fig 3). $\diamond$

\vskip0,5cm
 \paragraph{5. The main result}
        
\begin{flushleft}
 \textbf{Theorem}
\end{flushleft}
\emph{In a genus zero surface $\Sigma$ every free homotopy class   of loops with cobracket zero is a multiple of a simple class.}

\vskip0,5cm
\noindent\textbf{Proof.}
Hence $\Sigma= S^2\setminus   \left\{p_1, \ldots, p_k\  \right\}$  is a two sphere with $k $  punctures; clearly if $k<3$ every class is a multiple of a simple class.

Let $\Phi$ be a non simple class  of loops on $\Sigma$, and $\varphi : S^1 \longrightarrow \Sigma$ be a minimal generic representative of $\Phi$ .

Every $\varphi_{ab} $ which is a simple loop devides the sphere $S^2$ into two disks and so $\left\{p_1, \ldots, p_k\  \right\}$ into two proper parts. (  $\varphi$  cannot contain monogons since it is minimal ! )

Among them, choose $\varphi_{\overline{a}\overline{b}}$ so that one of the two parts in which $\left\{p_1, \ldots, p_k\  \right\}$ is partioned is a minimal subset. We can assume that this part is 
$\left\{p_1, \ldots, p_n\  \right\}$.  Put $\Psi = Class (\varphi_{\overline{a}\overline{b}})$

Since the cobracket of $\Phi$ is
$$\Delta (\Phi) = \Delta (\varphi) =\sum_{ a\approx b } sign (a,b)\cdot\  Class (\varphi_{ab}) \otimes Class (\varphi_{ba})$$

\noindent where  $sign (a,b)= 1$ if the couple of the positive tangent direction of $\varphi$ at $a$ and the positive tangent direction of $\varphi$ at $b$ is a positive pair according to the orientation of $\Sigma$ ; otherwise  $sign (a,b)= -1$. ( Remember that the definition of $sign (a,b)$ can be given also for each loop and that the cobracket is independent from the chosen representative of $\Phi$).

Since $\Delta (-\varphi) = -\Delta (\varphi)$, we can assume that   $\varphi_{\overline{\alpha}\overline{\beta}}$ goes around $\left\{p_1, \ldots, p_n\  \right\}$ counterclockwise.

So there is a couple $(a,b)$ with $a\approx b$ such that $Class\  \varphi_{ab}=\Psi$

We will show that, if $\Phi$ is  not a multiple of the simple class $\Psi$, then for each  couple $(a,b)$ with $a\approx b$ such that $Class\  \varphi_{ab}=\Psi$ is   $sign (a,b)= -1$.

\begin{center}

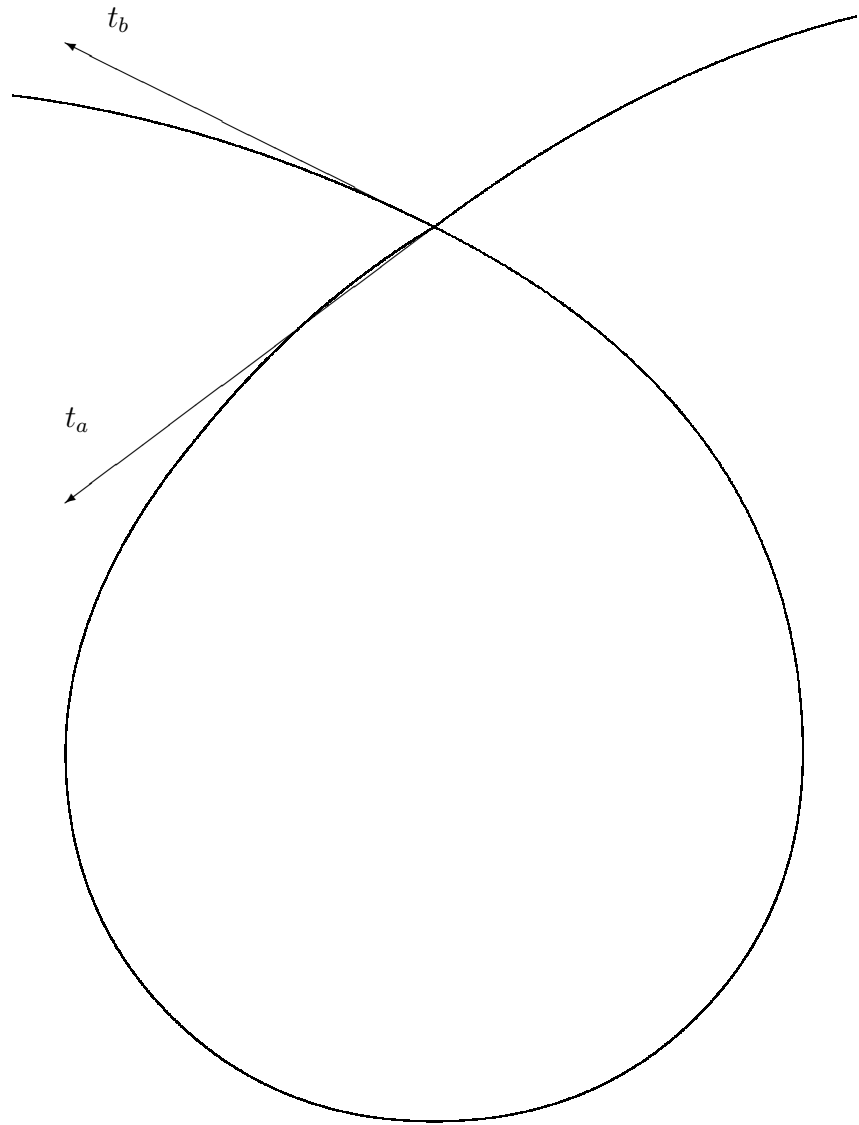
\textbf{Figure 7 }
\end{center}

\begin{center}

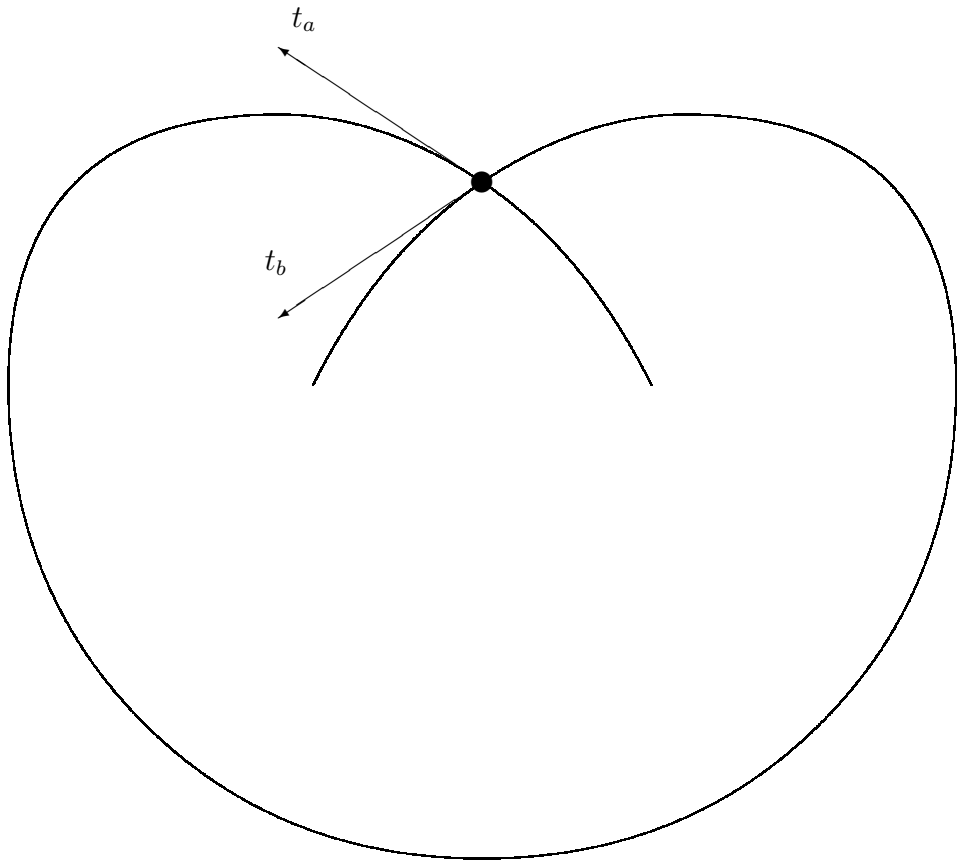
\textbf{Figure 8}
\end{center}

Since the free non-trivial homotopy classes of loops are free generators of the Lie bialgebra, this will show that $\Delta (\Phi)\neq 0$.

Hence, ab absurd, let us  suppose that $a\approx b$ with $Class\  \varphi_{ab}=\Psi$  and  $sign (a,b)= 1$.

Now, putting $p_k = \infty$, $\Sigma= \mathbb{R}^2\setminus   \left\{p_1, \ldots, p_{k-1}\  \right\}$ and  $\varphi_{ab}$ is a loop in  $\mathbb{R}^2$ going around $p_1, \ldots, p_n$  counterclockwise.
 
If all $im (\varphi)$ is  inside  $\varphi_{ab}$ then, choosing $\alpha$ and $\beta$ with $\alpha<b<a<\beta$ and  putting $f=\varphi_{\alpha\beta}$ we have that 
$f$ verifies the property of lemma 2 with $\psi=\varphi_{ab}$. 

Hence, if $f$ has  $m$ double points, in the notations of lemma, we have that $ a_m= b, b_m=a$ and  that $f_{ a_m b_m}=\varphi_{ba}\cong \psi^m$. Hence $\varphi\cong \psi^{m+1}$, contrarily to hypothesis.

Then suppose that not all  $im (\varphi)$ is  inside of $\varphi_{ab}$.
Since $sign (a,b)= 1$ we have the situation of figure 8, i.e., there is $c$ and $d $ so that $c<a<b<d$ with  $\varphi_{ca}$ and  $\varphi_{bc}$ inside of $\varphi_{ab}$. Choose the arc $cd $ maximal with this property. Then $c$ and $d$ are double points and $\varphi (c),\varphi (d)\in im (\varphi_{ab})$. 

Let $\psi$ be a simple loop, $\psi\cong\varphi_{ab}$ with $im (\varphi_{ab}) \subset D_ \psi$ so that not all  $im (\varphi)$ is  inside  $\psi$, and $im (\psi)$ contains no double point of $\varphi $. 

Take $s<c<d<t$ such that $\varphi (s),  \varphi (t)\in im (\psi)$ and $im (\varphi_{st})\subset D_ \psi$.

Hence $f=\varphi_{st}$ verifies the property of lemma 2 for $\psi$.

Again, in the notations of lemma,  we have that $ a_1= a, b_1=b$, and in $\varphi_{ab}=\varphi_{a_1b_1}$ there would be, by the lemma,  no double point of $f$ exept $f(a)=f(b)$, but $f(c)$ and $f(d)$ are. A contradiction. $\diamond$

\vskip1cm

\begin{flushleft}
 \textbf{ References}
\end{flushleft}

[1] Chas, Moira, Combinatorial Lie bialgebras of curves on surfaces, Topology 43 (2004) 543-568

[2] Goldman, William M., Invariant functions on Lie groups and Hamiltonian flows
of surface group representations. Invent. Math. 85 (1986), no. 2, 263-302

[3] Turaev, Vladimir G., Algebras of loops on surfaces, algebras of knots, and quantization, Braid Groups, Knot Theory and Statistical Mechanics, Adv.Ser.Math.Phys.,vol.9,World Sci.Publishing, Teaneck,NJ,1989, 59-95. 

[4] Turaev, Vladimir G., Skein quantization of Poisson algebras of loops on surfaces.
Ann. Sci. Ecole Norm. Sup.(4) 24(1991), no.6 , 635-704

\vfill
\break
\begin{flushleft}
 \textbf{ Figure Captions}
\end{flushleft}

\vskip1cm
fig. 1. Lemma 2, m=3
\vskip1cm
fig. 2.  m=1, $f(\beta)$ inside  $f_{a_1b_1}$
\vskip1cm
fig. 3. Case (1), m+1=3, $f(\beta)$  inside  $\sigma_m$
\vskip1cm
fig. 4. Case (2), m+1=4
\vskip1cm
fig. 5. Case (3), m+1=4
\vskip1cm
fig. 6. Case (4), m+1=4
\vskip1cm
fig. 7. $sign(a,b)=-1$
\vskip1cm
fig. 8. $sign(a,b)=1$
\vfill
\break

\setlength{\unitlength}{.4mm}

\begin{center}

 \begin{picture}(0,0)(0,150)
   \qbezier(140,140)(80,200)(0,200)
    \qbezier(200,0)(200,80)(140,140)
 \qbezier(200,0)(200,-80)(140,-140)
  \qbezier(140,-140)(80,-200)(0,-200)
    \qbezier(-200,0)(-200,80)(-140,140)
  \qbezier(-140,140)(-80,200)(0,200)
  \qbezier(-200,0)(-200,-80)(-140,-140)
  \qbezier(-140,-140)(-80,-200)(0,-200)

  \qbezier(140,140)(-50,100)(-50,0)
  \qbezier(-50,0)(-50,-50)(0,-50)
  \qbezier(50,0)(50,-50)(0,-50)
      \qbezier(-20,70)(50,70)(50,0)
    \qbezier(-20,70)(-80,70)(-80,0)
    \qbezier(-80,0)(-80,-80)(0,-80)
  \qbezier(80,0)(80,-80)(0,-80)
    \qbezier(-20,120)(80,120)(80,0)
  \qbezier(-20,120)(-140,120)(-140,0)
    \qbezier(-140,0)(-140,-130)(0,-130)
  \qbezier(130,0)(130,-130)(0,-130)
    \qbezier(-20,170)(130,170)(130,0)

  \put(0,0){\circle*{8}}\put(14,0){\circle*{8}}\put(5,14){\circle*{8}}
 \put(0,-280){fig. 1 }
\end{picture}
\end{center}

\vfill
\break

\setlength{\unitlength}{.4mm}

\begin{center}
\begin{picture}(0,0)(0,150)
   \qbezier(140,140)(80,200)(0,200)
   \qbezier(200,0)(200,80)(140,140)
    \qbezier(200,0)(200,-80)(140,-140)
  \qbezier(140,-140)(80,-200)(0,-200)
  
  \qbezier(-200,0)(-200,80)(-140,140)
 \qbezier(-140,140)(-80,200)(0,200)
  \qbezier(-200,0)(-200,-80)(-140,-140)
  \qbezier(-140,-140)(-80,-200)(0,-200)

  \qbezier(140,140)(-100,100)(-100,0)
     \qbezier(-100,0)(-100,-120)(0,-120)
  \qbezier(140,0)(140,-120)(0,-120)
    \qbezier(140,0)(140,60)(110,100)
  \qbezier(110,100)(88,129)(60,123)
    \qbezier(60,123)(32,117)(60,25)

  \put(0,0){\circle*{8}}\put(14,0){\circle*{8}}\put(5,14){\circle*{8}}\put(-8,15){\circle*{8}}

 \put(0,-290){fig. 2}
\end{picture}
\end{center}

\vfill
\break

\setlength{\unitlength}{.4mm}

\begin{center}
\begin{picture}(0,0)(0,150)
  \qbezier(200,0)(200,80)(140,140)
  \qbezier(140,140)(80,200)(0,200)
  \qbezier(200,0)(200,-80)(140,-140)
  \qbezier(140,-140)(80,-200)(0,-200)
  
  \qbezier(-200,0)(-200,80)(-140,140)
  \qbezier(-140,140)(-80,200)(0,200)
  \qbezier(-200,0)(-200,-80)(-140,-140)
  \qbezier(-140,-140)(-80,-200)(0,-200)

  \qbezier(140,140)(-50,100)(-50,0)
  \qbezier(-50,0)(-50,-50)(0,-50)
  \qbezier(50,0)(50,-50)(0,-50)
  
    \qbezier(-20,70)(50,70)(50,0)
    \qbezier(-20,70)(-80,70)(-80,0)
 
 
  \qbezier(-80,0)(-80,-80)(0,-80)
  \qbezier(80,0)(80,-80)(0,-80)
  
  \qbezier(-20,120)(80,120)(80,0)
  \qbezier(-20,120)(-140,120)(-140,0)

  \qbezier(-140,0)(-140,-130)(0,-130)
  \qbezier(140,0)(140,-130)(0,-130)
  
  \qbezier(140,0)(140,60)(110,100)
  \qbezier(110,100)(89,128)(75,122)
  \qbezier(75,122)(65,120)(95,70)

  \put(0,0){\circle*{8}}\put(14,0){\circle*{8}}\put(5,14){\circle*{8}}\put(-8,15){\circle*{8}}

 \put(0,-280){fig. 3}
\end{picture}
\end{center}

\vfill
\break


\begin{center}
\setlength{\unitlength}{.4mm}

\begin{picture}(0,0)(0,150)
 
  \qbezier(200,0)(200,80)(140,140)
  \qbezier(140,140)(80,200)(0,200)
  \qbezier(200,0)(200,-80)(140,-140)
  \qbezier(140,-140)(80,-200)(0,-200)
  
  \qbezier(-200,0)(-200,80)(-140,140)
  \qbezier(-140,140)(-80,200)(0,200)
  \qbezier(-200,0)(-200,-80)(-140,-140)
  \qbezier(-140,-140)(-80,-200)(0,-200)

  \qbezier(140,140)(-50,100)(-50,0)
  \qbezier(-50,0)(-50,-50)(0,-50)
  \qbezier(50,0)(50,-50)(0,-50)
      \qbezier(-20,70)(50,70)(50,0)
    \qbezier(-20,70)(-80,70)(-80,0)
    \qbezier(-80,0)(-80,-80)(0,-80)
  \qbezier(80,0)(80,-80)(0,-80)
    \qbezier(-20,120)(80,120)(80,0)
  \qbezier(-20,120)(-140,120)(-140,0)
    \qbezier(-140,0)(-140,-130)(0,-130)
  \qbezier(130,0)(130,-130)(0,-130)
    \qbezier(130,0)(130,150)(60,140)
    \qbezier(60,140)(25,135)(95,70)

  \put(0,0){\circle*{8}}\put(14,0){\circle*{8}}\put(5,14){\circle*{8}}
 \put(0,-280){fig. 4}
\end{picture}
\end{center}

\vfill
\break

\begin{center}

\setlength{\unitlength}{.4mm}
\begin{picture}(0,0)(0,150) 
  
  \qbezier(200,0)(200,80)(140,140)
  \qbezier(140,140)(80,200)(0,200)
  \qbezier(200,0)(200,-80)(140,-140)
  \qbezier(140,-140)(80,-200)(0,-200)
  
  \qbezier(-200,0)(-200,80)(-140,140)
  \qbezier(-140,140)(-80,200)(0,200)
  \qbezier(-200,0)(-200,-80)(-140,-140)
  \qbezier(-140,-140)(-80,-200)(0,-200)

  \qbezier(140,140)(-50,100)(-50,0)
  \qbezier(-50,0)(-50,-50)(0,-50)
  \qbezier(50,0)(50,-50)(0,-50)
  
    \qbezier(-20,70)(50,70)(50,0)
    \qbezier(-20,70)(-80,70)(-80,0)

  \qbezier(-80,0)(-80,-80)(0,-80)
  \qbezier(80,0)(80,-80)(0,-80)
  
  \qbezier(-20,120)(80,120)(80,0)
  \qbezier(-20,120)(-140,120)(-140,0)

  \qbezier(-140,0)(-140,-130)(0,-130)
  \qbezier(130,0)(130,-130)(0,-130)
  
  \qbezier(-20,170)(130,170)(130,0)
  \qbezier (-20,170)(-100,170)(-80,80)

  \put(0,0){\circle*{8}}\put(14,0){\circle*{8}}\put(5,14){\circle*{8}}

 \put(0,-280){fig. 5}
\end{picture}
\end{center}

\vfill
\break


\begin{center}
\setlength{\unitlength}{.4mm}
\begin{picture}(0,0)(0,150)

  \qbezier(200,0)(200,80)(140,140)
  \qbezier(140,140)(80,200)(0,200)
  \qbezier(200,0)(200,-80)(140,-140)
  \qbezier(140,-140)(80,-200)(0,-200)
  
  \qbezier(-200,0)(-200,80)(-140,140)
  \qbezier(-140,140)(-80,200)(0,200)
  \qbezier(-200,0)(-200,-80)(-140,-140)
  \qbezier(-140,-140)(-80,-200)(0,-200)

  \qbezier(140,140)(-50,100)(-50,0)
  \qbezier(-50,0)(-50,-50)(0,-50)
  \qbezier(50,0)(50,-50)(0,-50)
      \qbezier(-20,70)(50,70)(50,0)
    \qbezier(-20,70)(-80,70)(-80,0)
   \qbezier(-80,0)(-80,-80)(0,-80)
  \qbezier(80,0)(80,-80)(0,-80)
    \qbezier(-20,120)(80,120)(80,0)
  \qbezier(-20,120)(-140,120)(-140,0)
    \qbezier(-140,0)(-140,-130)(0,-130)
  \qbezier(130,0)(130,-130)(0,-130)
    \qbezier(-20,170)(130,170)(130,0)
  \qbezier(-20,170)(-40,170)(-40,150)
  \qbezier(-40,150)(-40,100)(40,180)
  
  \put(0,0){\circle*{8}}\put(14,0){\circle*{8}}\put(5,14){\circle*{8}}
 \put(0,-280){fig. 6}
\end{picture}
\end{center}

\vfill
\break

\begin{center}
\setlength{\unitlength}{.7mm}
\begin{picture}(0,0)(0,30)

\qbezier(0,0)(40,30)(80,40)
\qbezier(0,0)(-40,20)(-80,25)
\qbezier(0,0)(-25,-15)(-45,-40)
\qbezier(-45,-40)(-70,-70)(-70,-100)
\qbezier(70,-100)(70,-130)(50,-150)
\qbezier(50,-150)(30,-170)(0,-170)
\qbezier(0,0)(70,-35)(70,-100)
\qbezier(-70,-100)(-70,-130)(-50,-150)   
\qbezier(-50,-150)(-30,-170)(0,-170)
 \put(0,0){\vector(-2,1){70}}
 \put(0,0){\vector(-4,-3){70}} 
 \put(-70,-38){$t_a$}
 \put(-62,38){$t_b$}
\put(0,-210){fig. 7}
 \end{picture}  
\end{center}

\vfill
\break

\begin{center}
\setlength{\unitlength}{.9mm}
\begin{picture}(0,0)(0,30)

\qbezier(0,0)(15,10)(30,10)
\qbezier(30,10)(50,10)(60,0)
\qbezier(0,0)(-15,10)(-30,10)
\qbezier(-30,10)(-50,10)(-60,0)
\qbezier(60,0)(70,-10)(70,-30)
\qbezier(-60,0)(-70,-10)(-70,-30)
\qbezier(70,-30)(70,-60)(50,-80)
\qbezier(50,-80)(30,-100)(0,-100)
\qbezier(-70,-30)(-70,-60)(-50,-80)
\qbezier(-50,-80)(-30,-100)(0,-100)
\qbezier(0,0)(-15,-10)(-25,-30)
\qbezier(0,0)(15,-10)(25,-30)
 \put(0,0){\vector(-3,2){30}}
 \put(0,0){\vector(-3,-2){30}} 
 \put(-32,-13){$t_b$}
 \put(-28,23){$t_a$}
 \put(0,0){\circle*{3}}
\put(0,-145){fig. 8}

\end{picture}
\end{center}

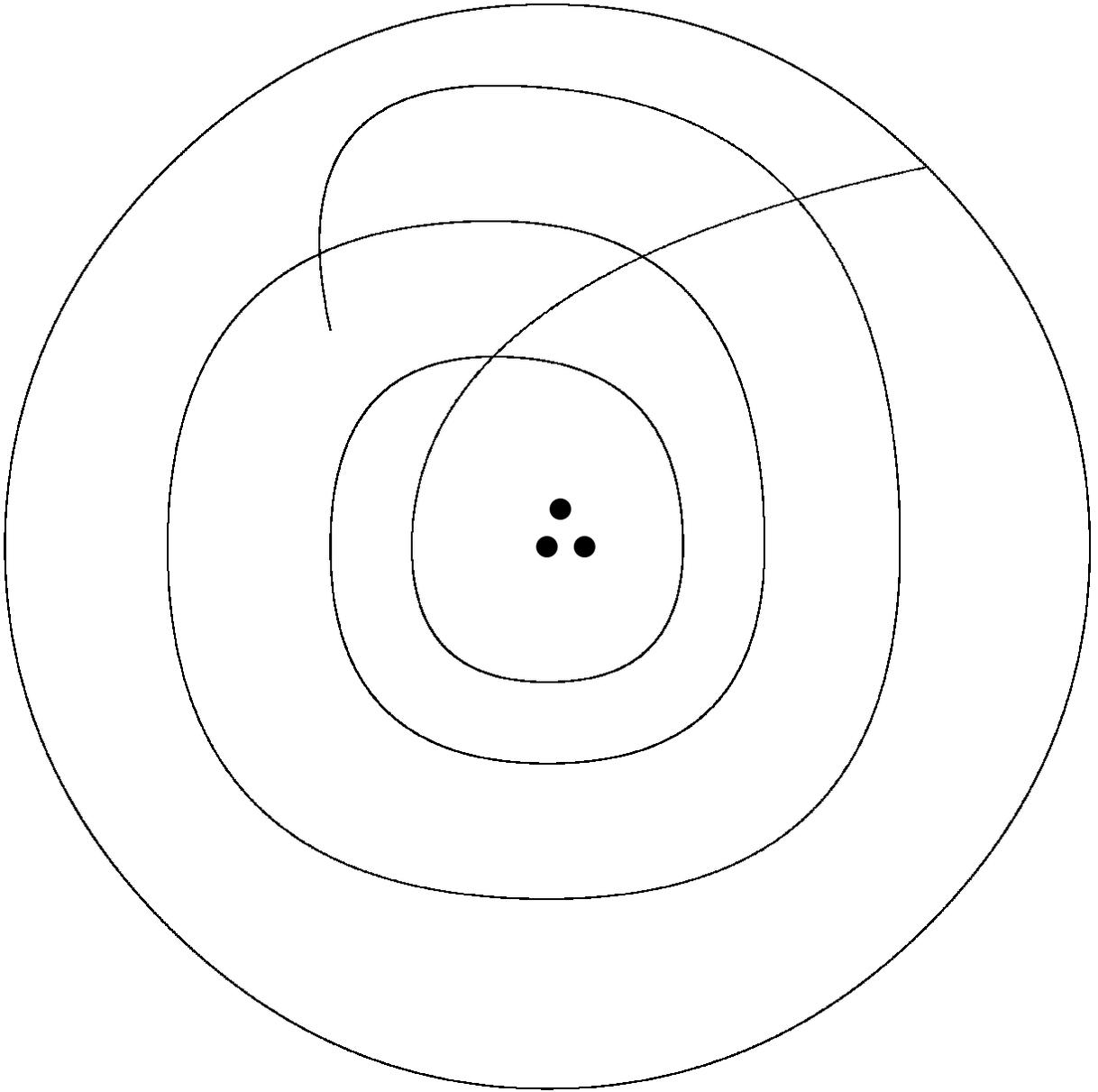
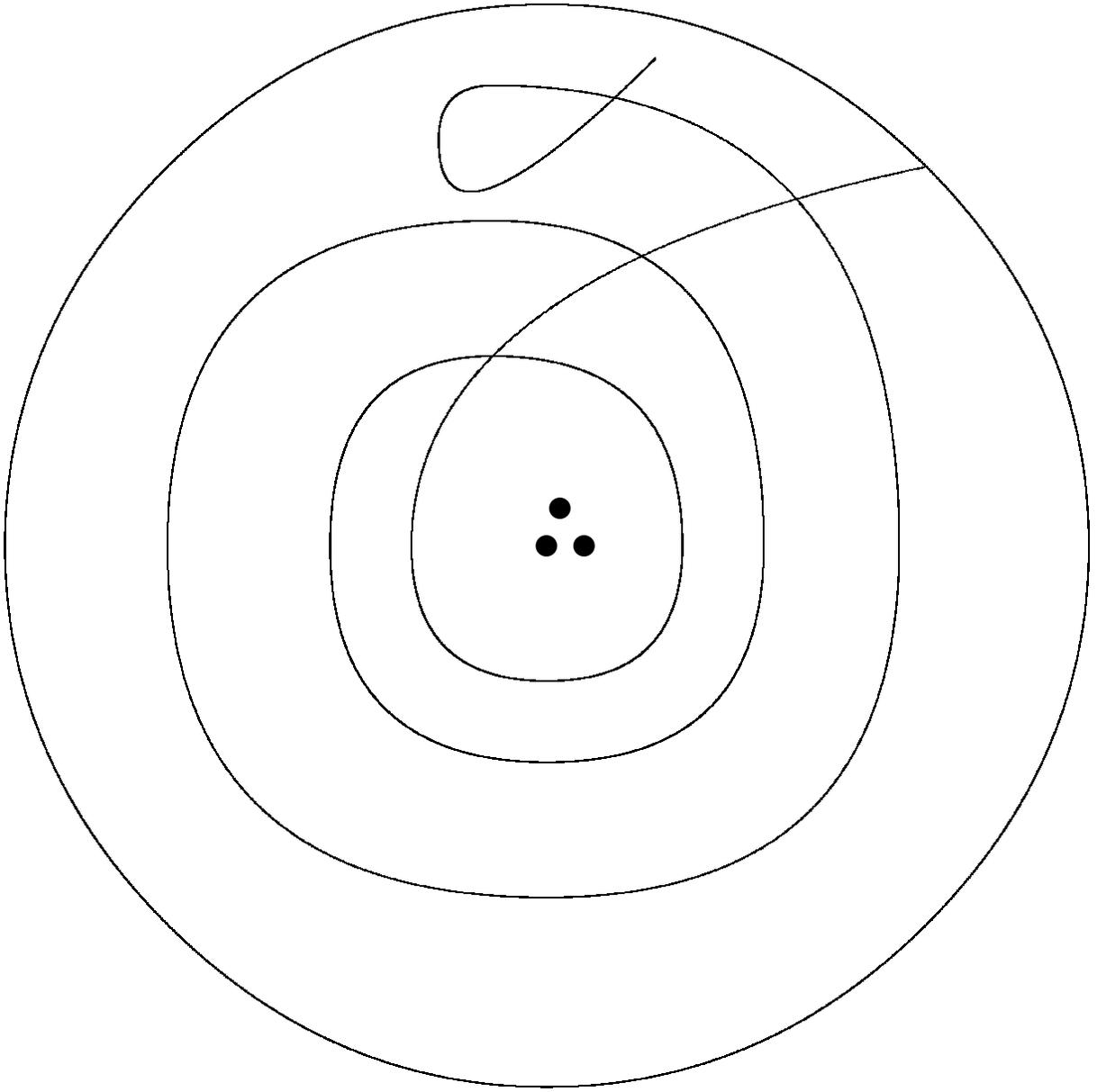
\end{document}